\documentclass[12pt,reqno,a4wide]{amsart}

\usepackage{tikz} 
\usepackage{calrsfs}
\usepackage{mathrsfs}

\usepackage{array}
\usepackage{hyperref}

\usetikzlibrary{decorations.pathreplacing}
\usetikzlibrary{fit}

\usepackage{pgfplots}

\allowdisplaybreaks

\begin{document}
\bibliographystyle{plain}

%
%

	\title[Prefixes of Stanley's Catalan paths with odd returns to the $x$-axis]
	{Prefixes of Stanley's Catalan paths with odd returns to the $x$-axis -- standard version and skew Catalan-Stanley paths}

	\author[H. Prodinger ]{Helmut Prodinger }
	\address{Department of Mathematics, University of Stellenbosch 7602, Stellenbosch, South Africa
	and
NITheCS (National Institute for
Theoretical and Computational Sciences), South Africa.}
	\email{hproding@sun.ac.za}

	\keywords{Catalan (Dyck) paths, odd returns to $x$-axis, skew paths, prefix, kernel method}
	
	\begin{abstract}
	Stanley considered Dyck paths where each maximal run of down-steps to the $x$-axis has odd length; they are also enumerated by (shifted)
Catalan numbers. Prefixes of these combinatorial objects are enumerated using the kernel method. A more challenging version of skew Dyck paths 
combined with Stanley's restriction is also considered.
	\end{abstract}
	
	\subjclass[2010]{05A15}

\maketitle

\section{Introduction}

Dyck paths consist of up-steps $(1,1)$ and down-steps $(1,-1)$, are not allowed to go below the $x$-axis and end at the $x$-axis.
The enumeration (path consisting of $2n$ steps) is via the Catalan numbers $\mathcal{C}_n=\frac1{n+1}\binom{2n}{n}$. Sometimes, paths that do not necessarily
end at the $x$-axis are also called Dyck (Catalan) paths. The ultimate reference is Stanley's book \cite{Stanley-book}. This book offers in Exercise 26 an
interesting alternative: Maximal sequences of down-steps, when the end on the $x$-axis must always have odd length. This restriction leads to smaller numbers of 
paths, but interestingly restricted paths consisting of $2n+2$ steps are also enumerated by Catalan numbers $\mathcal{C}_n$. The 5 allowed paths of length 8 are depicted in Figure~\ref{fig-5}:
	\begin{figure}[h]\label{fig-5}
	\begin{center}
		\begin{tikzpicture}[scale=0.5]
			\foreach \x in {0,1,2,3,4,5,6,7,8}
			{
				\draw[thin ] (\x,0) to  (\x,3);	
			}
			\foreach \x in {0,1,2,3}
			{
				\draw[thin ] (0,\x) to  (8,\x);	
			}		
			\draw[ultra thick](0,0) to (1,1);
			\draw[ultra thick](2,0) to (3,1);
			\draw[ultra thick](4,0) to (5,1);
			\draw[ultra thick](6,0) to (7,1);
			
			\draw[ultra thick,red] (1,1) to (2,0);
			\draw[ultra thick,red] (3,1) to (4,0);
			\draw[ultra thick,red] (5,1) to (6,0);
			\draw[ultra thick,red] (7,1) to (8,0);
			
		\end{tikzpicture}\ %
		\begin{tikzpicture}[scale=0.5]
			\foreach \x in {0,1,2,3,4,5,6,7,8}
			{
				\draw[thin ] (\x,0) to  (\x,3);	
			}
			\foreach \x in {0,1,2,3}
			{
				\draw[thin ] (0,\x) to  (8,\x);	
			}		
			\draw[ultra thick](0,0) to (1,1);
			\draw[ultra thick](2,0) to (5,3);
			
			\draw[ultra thick,red] (1,1) to (2,0);
			\draw[ultra thick,red] (5,3) to (8,0);
			
		\end{tikzpicture}\ %
		\begin{tikzpicture}[scale=0.5]
			\foreach \x in {0,1,2,3,4,5,6,7,8}
			{
				\draw[thin ] (\x,0) to  (\x,3);	
			}
			\foreach \x in {0,1,2,3}
			{
				\draw[thin ] (0,\x) to  (8,\x);	
			}		
			\draw[ultra thick](0,0) to (2,2);
			\draw[ultra thick](3,1) to (5,3);
			
			\draw[ultra thick] (2,2) to (3,1);
			\draw[ultra thick,red] (5,3) to (8,0);
					\end{tikzpicture} 
				\vskip0.1cm
		\begin{tikzpicture}[scale=0.5]
			\foreach \x in {0,1,2,3,4,5,6,7,8}
			{
				\draw[thin ] (\x,0) to  (\x,3);	
			}
			\foreach \x in {0,1,2,3}
			{
				\draw[thin ] (0,\x) to  (8,\x);	
			}		
			\draw[ultra thick](0,0) to (3,3);
			\draw[ultra thick](6,0) to (7,1);
			
			\draw[ultra thick,red] (3,3) to (6,0);
			\draw[ultra thick,red] (7,1) to (8,0);
			
		\end{tikzpicture}\ %
		\begin{tikzpicture}[scale=0.5]
			\foreach \x in {0,1,2,3,4,5,6,7,8}
			{
				\draw[thin ] (\x,0) to  (\x,3);	
			}
			\foreach \x in {0,1,2,3}
			{
				\draw[thin ] (0,\x) to  (8,\x);	
			}		
			\draw[ultra thick](0,0) to (3,3);
			\draw[ultra thick](4,2) to (5,3);
			
			\draw[ultra thick] (3,3) to (4,2);
			\draw[ultra thick,red] (5,3) to (8,0);
			
		\end{tikzpicture}
	\end{center}
\caption{Restricted paths \`a la Stanley of length 8. The maximal runs of down-steps to the $x$-axis are depicted in red.}
\end{figure}
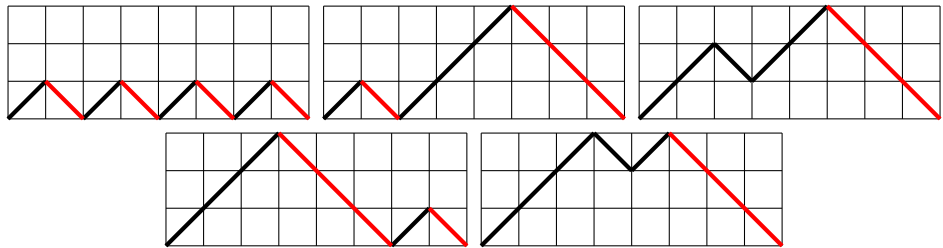

 In the first of this paper we will enumerate prefixes of restricted paths \`a la Stanley: they can be completed to allowed paths.
 An alternative name of these combinatorial objects is \emph{partial} Stanley-Dyck paths.
 We might also say that a maximal down-run  the $x$-axis of \emph{even} length has never occured in a partial Dyck path. The method of choice is
 generating functions (in $z$ marking the length) but also in $u$, marking the final level. The details will appear in the next section.
 
 In the final section an additional restriction (skew Dyck paths) will be combined with the Stanley restriction: They allow an
 additional type of south-west step, but no overlap is allowed to occur. As in a previous publication  \cite{Prodinger-prefix} we prefer to
 draw a down-step $(1,-1)$ but depict it with the red color.
 
 As in various examples discussed in the past \cite{Prodinger-kernel, Prodinger-prefix, garden} we will use the kernel method to set up appropriate
 generating functions. There are three bivariate generating functions $F(u), G(u), H(u)$ (also depending on $z$) related to the nature of the last step of the 
 prefix of a Stanley-Dyck path. To solve the system, a `bad' factor $u-r_2$ must be divided out from numerator and denominator, after which one can plug in $u=0$
 and identify the unknown quantities $f_1=F(0)$, $h_0=H(0)$. The function $r_2(z)$ generates Catalan numbers.
 
 In the final section, instead of 3 functions, we have to deal with 5 functions, but the strategy is similar. The function that plays the role of $r_2$ now has coefficients
 1, 3, 10, 36, 137, 543, 2219, 9285, 39587, 171369, 751236, which is sequence A2212 in \cite{OEIS}. The paper \cite{garden} contains a fair amount of combinatorics around this sequence.
 
\section{Prefixes of Dyck paths with only maximal down-runs to the $x$-axis of odd length} 
 
 We present a graph that allows to identify all prefixes of Stanley-Dyck path. It has 3 layers of states, and maximal sequences of down-steps
 are depicted in blue. Parity plays a role: if such a sequence would lead to the $x$-axis after an even number of states, this is prohibited as
 such a state in the first layer is missing. Otherwise we introduce sequences $f_i$, $i\ge1$,  $g_i$, $h_i$, $i\ge0$. These quantities all depend on the variable $z$ and
 describe generating functions of paths leading to a particular state. The state $h_0$ is special, it is related to Stanley-Dyck paths.
  \begin{figure}[h]

 	\begin{center}
 		\begin{tikzpicture}[scale=1.8,main node/.style={circle,draw,font=\Large\bfseries}]

 			\foreach \x in {0,1,2,3,4,5,6,7,8}
 			{
 				\draw (\x,0) circle (0.05cm);
 				\fill (\x,0) circle (0.05cm);
 				\draw (\x,-1) circle (0.05cm);
 				\fill (\x,-1) circle (0.05cm);
 			}
 			
 			\foreach \x in {1,2,3,4,5,6,7,8}
{
	\draw (\x,0) circle (0.05cm);
	\fill (\x,0) circle (0.05cm);
	\draw (\x,1) circle (0.05cm);
	\fill (\x,1) circle (0.05cm);
}

 			\fill (0,0) circle (0.08cm);

 			\foreach \x in {0,2,4,6}
 			{
 			}

 			\foreach \x in {0,1,2,3,4,5,6,7}
 			{
 				\draw[ultra thick, -latex] (\x,0) to  (\x+1,0);	
\draw[ultra thick,cyan, latex-] (\x,-1) to  (\x+1,-1);	

 			}			
 			
 			\foreach \x in {1,2,3,4,5,6,7}
 			{
 				\draw[ultra thick,cyan, latex-] (\x,1) to  (\x+1,1);	
 				
 			}

 			 			\foreach \x in {0,2,4,6}
 			 			{ 				\draw[ultra thick,cyan, -latex] (\x+1,0) to  (\x,-1);	
 			 			}
 		 			
 		 						\foreach \x in {0,2,4,6}
 		 			{ 				\draw[ultra thick,cyan, -latex] (\x+2,0) to  (\x+1,1);	
 		 			}

 			\foreach \x in {0,1,2,3,4,5,6,7}
 			{ 				
 				\draw[ultra thick, latex-] (\x+1,0)[out=200, in=80] to  (\x,-1);	
 			}
 			
 				\foreach \x in {1,2,3,4,5,6,7}
 			{ 				
 				\draw[ultra thick, latex-] (\x+1,0)[in=-20, out=100] to  (\x,1);	
 			}

 			\foreach \x in {0,1,2,3,4,5,6,7,8}
 			{
 			}


 			\foreach \x in {0,1,2,3,4,5,6,7,8}
{
	\draw (\x,0) circle (0.05cm);
	\fill (\x,0) circle (0.05cm);
	\draw (\x,-1) circle (0.05cm);
	\fill (\x,-1) circle (0.05cm);
}

\foreach \x in {1,2,3,4,5,6,7,8}
{
	\draw (\x,0) circle (0.05cm);
	\fill (\x,0) circle (0.05cm);
	\draw (\x,1) circle (0.05cm);
	\fill (\x,1) circle (0.05cm);
}

\fill (0,0) circle (0.08cm);

 		\end{tikzpicture}
 	\caption{Three layers of states, labelled $f,g,h$, in that order. The state $h_0$ is responsible for Stanley-Dyck paths, and all others to prefixes of them.}
 	\end{center}
 \end{figure}
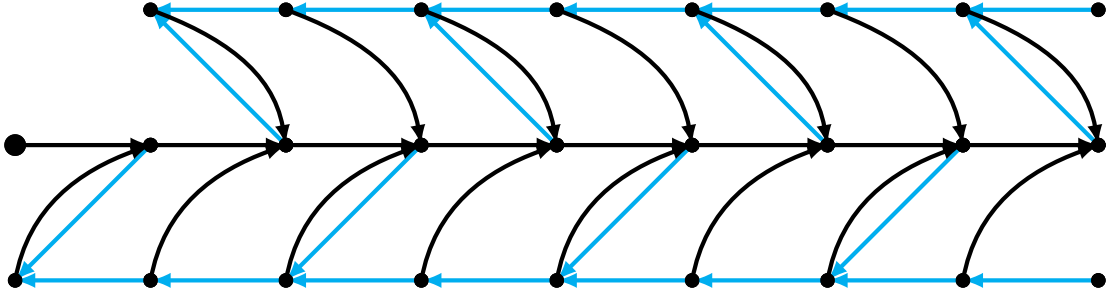
We read off recursions and introduce bivariate generating functions
\begin{equation*}
F(u)=\sum_{i\ge1}f_iu^{i-1},\quad G(u)=\sum_{i\ge0}g_iu^{i},\quad H(u)=\sum_{i\ge0}h_iu^{i}.
\end{equation*}
The recursions are according to the last symbol read. Firstly,
\begin{align*}
f_i&=zf_{i+1}+[i\text{ odd}]zg_{i+1},\ i\ge1,\quad\text{and}\\
	\sum_{i\ge1}f_iu^i&=\sum_{i\ge1}zf_{i+1}u^i+\sum_{i\ge1}[i\text{ odd}]zg_{i+1}u^i, \\
	uF(u)&=zF(u)-zf_1+\frac zu\sum_{i\ge0}[i\text{ even}]g_{i}u^i-\frac zu, \\
	uF(u)&=zF(u)-zf_1+\frac z{2u}G(u)+\frac z{2u}G(-u)-\frac zu.
\end{align*}
Secondly,
\begin{align*}
	h_i&=zh_{i+1}+[i\text{ even}]zg_{i+1},\ i\ge0,\quad\text{and}\\
H(u)&=\frac zuH(u)-\frac zuh_0+\frac zu\sum_{i\ge1}[i\text{ odd}]g_{i}u^{i},\\
	H(u)&=\frac zuH(u)-\frac zuh_0+\frac z{2u}G(u)-\frac z{2u}G(-u).
\end{align*}
Thirdly
\begin{align*}
g_{i+1}&=zf_i+zg_i+zh_i,\ i\ge1,\ g_1=z+zh_0,\ g_0=1,\quad\text{and}\\
	\sum_{i\ge1}g_{i+1}u^i&=\sum_{i\ge1}zf_iu^i+\sum_{i\ge1}zg_iu^i+\sum_{i\ge1}zh_iu^i, \ g_1=z+zh_0,\ g_0=1,\\
	\frac1u G(u)-\frac1u&=zuF(u)+zG(u)+zH(u).
\end{align*}
The 3 equations can be solved: We start with the simplest function
\begin{equation*}
G(u)=\frac{z^2u^2f_1+z^2uh_0-u+z^2u+z}{z-u+zu^2};
\end{equation*}
to get a form that leads to further developments, we factorize the denominator:
\begin{equation*}
z-u+zu^2=z(u-r_1)(u-r_2),\quad r_1=\frac{1+\sqrt{1-4z^2}}{2z},\quad r_2=\frac{1-\sqrt{1-4z^2}}{2z}.
\end{equation*}
Dividing out the factor $u-r_2$  from numerator and denominator, we find\footnote{This is in accordance with the kernel method.}
\begin{equation*}
G(u)=\frac{r_2z^2f_1+z^2h_0-1+z^2+z^2uf_1}{r_2z-1+zu},\quad\text{and}\quad G(0)=1=\frac{r_2z^2f_1+z^2h_0-1+z^2}{r_2z-1}.
\end{equation*}
Using this for $G(u)$ and $G(-u)$, we get
\begin{equation*}
F(u)=\frac{\Xi}{2u(z-u+zu^2)(u-z)}
\end{equation*}
with $\Xi=z(-2zf_1u^3+z^2u^2f_1+zu^2G(-u)-2zu^2+2u^2f_1+z^2uh_0+z^2u-2zf_1u-uG(-u)+u+zG(-u)-z)$.
Dividing out the factor $u(u-z)(u-r_2)$ we find an attractive answer
\begin{equation*}
F(u)=\frac{z^3(1+uf_1)}{1-z^2-zr_2-z^2u^2}, \quad F(0)=f_1=\frac{z^3}{1-z^2-zr_2}=\frac{1-2z^2-\sqrt{1-4z^2}}{2z}=r_2-z.
\end{equation*}
From the equation $G(0)=1=\dfrac{r_2z^2f_1+z^2h_0-1+z^2}{r_2z-1}$, with the already known $f_1$, we conclude
\begin{equation*}
h_0=\frac{1-\sqrt{1-4z^2}}{2}=zr_2=z^2\sum_{n\ge0} \mathcal{C}_nz^{2n}=\sum_{n\ge0} \mathcal{C}_nz^{2n+2},
\end{equation*}
as predicted by Stanley \cite{Stanley-book}. Further simplifications include
\begin{equation*}
	F(u)=\frac{z(\frac{r_2}{z}-1)(1+zu(\frac{r_2}{z}-1))}{1-(\frac{r_2}{z}-1)u^2}=\frac{(1+uzr_2^2)r_2^2}{1-r_2^2u^2},
\end{equation*}
\begin{align*}
	H(u)&=\frac{z^2(1-z^2+uzh_0-zf_1)}{1-z^2-r_2z-z^2u^2}\\
	&=\frac{z^2(1-z^2+uz^2r_2-z(r_2-z))}{1-z^2-r_2z}\frac{1}{1-(\frac{r_2}{z}-1)u^2}=\frac{(1+uzr_2^2)zr_2}{1-r_2^2u^2}
\end{align*}
and
\begin{equation*}
	G(u)=\frac{1-uz^2r_2^3}{1-r_2u}.
\end{equation*}
This allows the following expansions
\begin{gather*}
	F(u)=\frac{(1+uzr_2^2)r_2^2}{1-r_2^2u^2}
	=\sum_{k\ge0}r_2^{2k+2}u^{2k}+z\sum_{k\ge0}r_2^{2k+4}u^{2k+1},\\
	G(u)=\frac{1-uz^2r_2^3}{1-r_2u}=\sum_{k\ge0}r_2^ku^k-z^2\sum_{k\ge0}r_2^{k+3}u^{k+1},\\
	H(u)=\frac{(1+uzr_2^2)zr_2}{1-r_2^2u^2}=z\sum_{k\ge0}r_2^{2k+1}u^{2k}+z^2\sum_{k\ge0}r_2^{2k+3}u^{2k+1}.
\end{gather*}
Classically, we know that
\begin{equation*}
r_2^j=z^j\sum_{n\ge0}\frac{j(2n+j-1)!}{n!(n+j)!}z^{2n}.
\end{equation*}
Then
\begin{equation*}
f_j=[u^{j-1}]F(u)=\begin{cases}r_2^{j+1}& j\text{ odd},\\
	zr_2^{j+2}& j\text{ even}.
	\end{cases}
\end{equation*}
Further, $g_0=1$ and for $j\ge1$
\begin{equation*}
g_j=r_2^j-z^2r_2^{j+2}.
\end{equation*}
Finally,
\begin{equation*}
h_j=[u^j]z\sum_{k\ge0}r_2^{2k+1}u^{2k}+[u^j]z^2\sum_{k\ge0}r_2^{2k+3}u^{2k+1}
=\begin{cases}z^2r_2^{j+2}& j\text{ odd},\\
	zr_2^{j+1}& j\text{ even}.
\end{cases}
\end{equation*}

An interesting quantity is $F(1)+G(1)+H(1)$, which is the generating function of \emph{all} prefixes of Stanley-Dyck paths, regardless of where they end;
\begin{align*}
F(1)&+G(1)+H(1)=\frac{1-z}{2}+\frac12\frac{1+z}{1-2z}\sqrt{1-4z^2}\\
&=1+z+2z^2+3z^3+5z^4+9z^5+16z^6+30z^7+55z^8+105z^9+196z^{10}+378z^{11}+\dots\,.
\end{align*}

\section{Skew Dyck paths and the Stanley restriction}

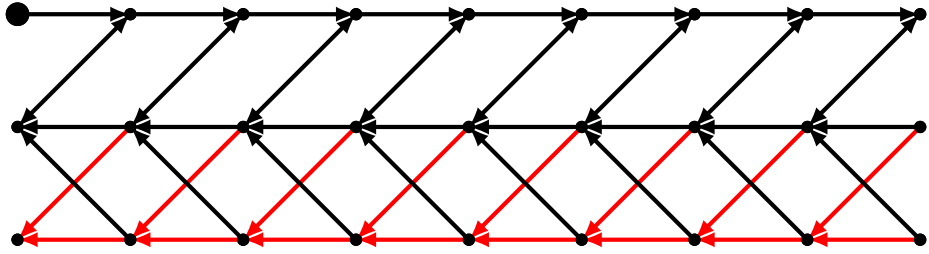
\begin{figure}[h]\label{figskew}

	\begin{center}
		\begin{tikzpicture}[scale=1.5]
			\draw (0,0) circle (0.1cm);
			\fill (0,0) circle (0.1cm);
			
			\foreach \x in {0,1,2,3,4,5,6,7,8}
			{
				\draw (\x,0) circle (0.05cm);
				\fill (\x,0) circle (0.05cm);
			}
			
			\foreach \x in {0,1,2,3,4,5,6,7,8}
			{
				\draw (\x,-1) circle (0.05cm);
				\fill (\x,-1) circle (0.05cm);
			}
			
			\foreach \x in {0,1,2,3,4,5,6,7,8}
			{
				\draw (\x,-2) circle (0.05cm);
				\fill (\x,-2) circle (0.05cm);
			}
			
			\foreach \x in {0,1,2,3,4,5,6,7}
			{
				\draw[ ultra thick,-latex] (\x,0) -- (\x+1,0);
				
			}

			\foreach \x in {0,1,2,3,4,5,6,7}
			{
				\draw[ ultra thick,latex-latex] (\x+1,0) -- (\x,-1);
			}
			
			\foreach \x in {0,1,2,3,4,5,6,7}
			{
				\draw[ ultra thick,-latex] (\x+1,-1) -- (\x,-1);
				
			}
			\foreach \x in {0,1,2,3,4,5,6,7}
			{
				\draw[ ultra thick,-latex,red] (\x+1,-1) -- (\x,-2);
				
			}
			
			\foreach \x in {0,1,2,3,4,5,6,7}
			{
				\draw[ ultra thick,-latex,red] (\x+1,-2) -- (\x,-2);
				
			}
			
			\foreach \x in {0,1,2,3,4,5,6,7}
			{
				\draw[ ultra thick,-latex] (\x+1,-2) -- (\x,-1);
				
			}

			\draw (0,0) circle (0.1cm);
\fill (0,0) circle (0.1cm);

\foreach \x in {0,1,2,3,4,5,6,7,8}
{
	\draw (\x,0) circle (0.05cm);
	\fill (\x,0) circle (0.05cm);
}

\foreach \x in {0,1,2,3,4,5,6,7,8}
{
	\draw (\x,-1) circle (0.05cm);
	\fill (\x,-1) circle (0.05cm);
}

\foreach \x in {0,1,2,3,4,5,6,7,8}
{
	\draw (\x,-2) circle (0.05cm);
	\fill (\x,-2) circle (0.05cm);
}

		\end{tikzpicture}
	\caption{Three layers of states according to skew Dyck paths.}
	\end{center}

\end{figure}

Skew Dyck paths \cite{Emeric} are augmented Dyck paths with additional south-west steps $(-1,-1)$, but no overlaps. In \cite{Prodinger-prefix} we preferred to
replace these extra steps by red down-steps; red-up resp.\ up-red are forbidden because of overlaps. The graph Figure~\ref{figskew} describes the allowed paths. If
they end on level 0, they are skew Dyck paths, otherwise just prefixes of them.

\begin{figure}[h]

	\begin{center}
		\begin{tikzpicture}[scale=1.8,main node/.style={circle,draw,font=\Large\bfseries}]

			\foreach \x in {0,1,2,3,4,5,6,7,8}
			{
				\draw (\x,0) circle (0.05cm);
				\fill (\x,0) circle (0.05cm);
				\draw (\x,-1) circle (0.05cm);
				\fill (\x,-1) circle (0.05cm);
			}
			
			\foreach \x in {1,2,3,4,5,6,7,8}
			{
				\draw (\x,0) circle (0.05cm);
				\fill (\x,0) circle (0.05cm);
				\draw (\x,1) circle (0.05cm);
				\fill (\x,1) circle (0.05cm);
			}

			\fill (0,0) circle (0.08cm);

			\foreach \x in {0,2,4,6}
			{
			}

			\foreach \x in {0,1,2,3,4,5,6,7}
			{
				\draw[ultra thick, -latex] (\x,0) to  (\x+1,0);	
				\draw[ultra thick,cyan, latex-] (\x,-1) to  (\x+1,-1);	
				
			}			
			
			\foreach \x in {1,2,3,4,5,6,7}
			{
				\draw[ultra thick,cyan, latex-] (\x,1) to  (\x+1,1);	
				
			}

			\foreach \x in {0,2,4,6}
			{ 				\draw[ultra thick,cyan, -latex] (\x+1,0) to  (\x,-1);	
			}
			
			\foreach \x in {0,2,4,6}
			{ 				\draw[ultra thick,cyan, -latex] (\x+2,0) to  (\x+1,1);	
			}

			\foreach \x in {0,1,2,3,4,5,6,7}
			{ 				
				\draw[ultra thick, latex-] (\x+1,0)[out=200, in=80] to  (\x,-1);	
			}
			
			\foreach \x in {1,2,3,4,5,6,7}
			{ 				
				\draw[ultra thick, latex-] (\x+1,0)[in=-20, out=100] to  (\x,1);	
			}

			\foreach \x in {1,2,3,4,5,6,7,8}
			{
					\draw[ultra thick, -latex,red] (\x,-1) to  (\x-1,-2);	
										\draw[ultra thick, -latex,red] (\x,-2) to  (\x-1,-2);	
																				\draw[ultra thick, -latex,cyan] (\x,-2) to  (\x-1,-1);	
			}			
			\foreach \x in {2,3,4,5,6,7,8}
			{
				\draw[ultra thick, -latex,red] (\x,1) to  (\x-1,2);	
				\draw[ultra thick, -latex,red] (\x,2) to  (\x-1,2);	
				\draw[ultra thick, -latex,cyan] (\x,2) to  (\x-1,1);	
			}	
			
			
			\foreach \x in {0,1,2,3,4,5,6,7,8}
			{
				\draw (\x,0) circle (0.05cm);
				\fill (\x,0) circle (0.05cm);
				\draw (\x,-1) circle (0.05cm);
				\fill (\x,-1) circle (0.05cm);
				\draw (\x,-2) circle (0.05cm);
				\fill (\x,-2) circle (0.05cm);
			}
			
			\foreach \x in {1,2,3,4,5,6,7,8}
			{
				\draw (\x,0) circle (0.05cm);
				\fill (\x,0) circle (0.05cm);
				\draw (\x,1) circle (0.05cm);
				\fill (\x,1) circle (0.05cm);
				\draw (\x,2) circle (0.05cm);
\fill (\x,2) circle (0.05cm);

			}

			\fill (0,0) circle (0.08cm);

		\end{tikzpicture}
	\caption{Graph to recognize skew Stanley-Dyck paths.}
	\end{center}
\end{figure}
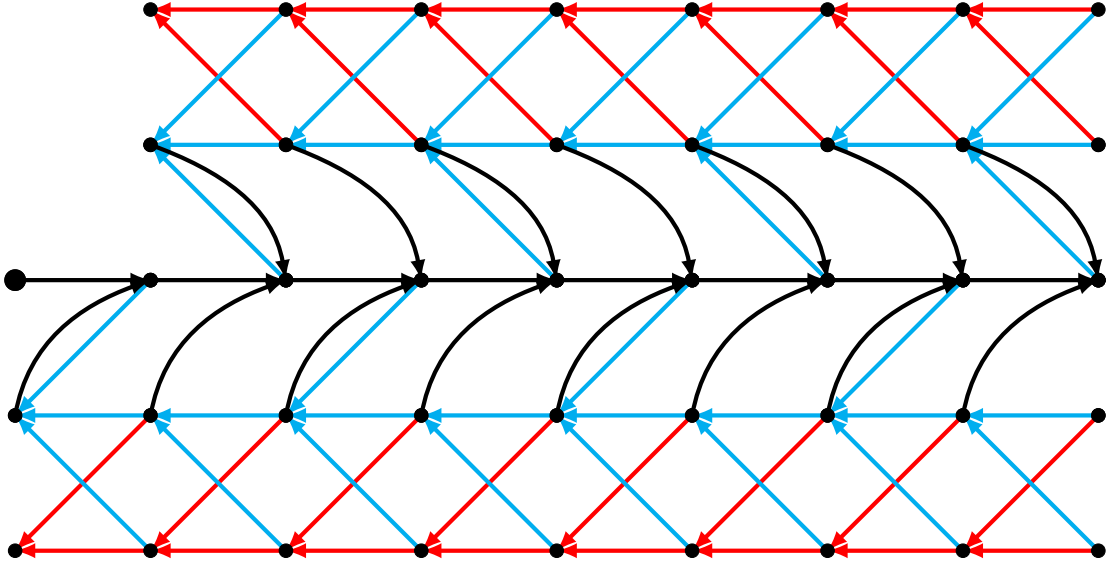
The relevant recursions can be read off, again by considering a last step of a partial path.
 \begin{align*}
e_i&=ze_{i+1}+zf_{i+1}, i\ge1\\
f_i&=ze_{i+1}+zf_{i+1}+[i\text{ odd}]zg_{i+1}, i\ge1\\
g_{i+1}&=zf_i+zg_i+zh_i,\ i\ge1,\ g_1=z+zh_0,\ g_0=1\\
h_i&=zh_{i+1}+[i\text{ even}]zg_{i+1}+zk_{i+1},\ i\ge0\\
k_i&=zh_{i+1}+zk_{i+1},\ i\ge0.
 \end{align*}
Summing, 
\begin{align*}
	\sum_{i\ge1}u^ie_i&=\sum_{i\ge1}u^ize_{i+1}+\sum_{i\ge1}u^izf_{i+1}\\
	\frac uzE(u)&=E(u)-e_1+F(u)-f_1;
\end{align*} 	
\begin{align*}
\sum_{i\ge1}u^if_i&=\sum_{i\ge1}u^ize_{i+1}+\sum_{i\ge1}u^izf_{i+1}+\sum_{i\ge1}u^i[i\text{ odd}]zg_{i+1}\\
\frac uzF(u)&=E(u)-e_1+F(u)-f_1+\frac1u\sum_{i\ge2}u^i[i\text{ even }]g_{i}\\
\frac uzF(u)&=E(u)-e_1+F(u)-f_1+\frac1{2u}G(u)+\frac1{2u}G(-u)-\frac1u;
\end{align*} 
\begin{align*}
g_1+\sum_{i\ge1}u^ig_{i+1}&=\sum_{i\ge1}u^izf_i+\sum_{i\ge1}u^izg_i+\sum_{i\ge1}u^izh_i+z+zh_0\\
\sum_{i\ge0}u^ig_{i+1}&=\sum_{i\ge1}u^izf_i+\sum_{i\ge0}u^izg_i+\sum_{i\ge0}u^izh_i\\
\frac{G(u)-1}{u}&=uzF(u)+zG(u)+zH(u);
\end{align*} 
\begin{align*}
\sum_{i\ge0}u^{i+1}h_i&=\sum_{i\ge0}u^{i+1}zh_{i+1}+\sum_{i\ge0}u^{i+1}[i\text{ even}]zg_{i+1}+\sum_{i\ge0}u^{i+1}zk_{i+1}\\
uH(u)&=z(H(u)-h_0)+z\sum_{i\ge1}u^{i}[i\text{ odd}]g_{i}+z(K(u)-k_0)\\
uH(u)&=z(H(u)-h_0)+\frac z2G(u)-\frac z2G(-u) +z(K(u)-k_0);
\end{align*}
\begin{align*}
	\sum_{i\ge0}u^{i+1}k_i&=\sum_{i\ge0}u^{i+1}zh_{i+1}+\sum_{i\ge0}u^{i+1}zk_{i+1}\\
	uK(u)&=z(H(u)-h_0)+z(K(u)-k_0).
\end{align*}
For convenience, we repeat the system that needs to be solved:
\begin{align*}
		\frac uzE(u)&=E(u)-e_1+F(u)-f_1,\\
		\frac uzF(u)&=E(u)-e_1+F(u)-f_1+\frac1{2u}G(u)+\frac1{2u}G(-u)-\frac1u,\\
		\frac{G(u)-1}{u}&=uzF(u)+zG(u)+zH(u),\\
		uH(u)&=z(H(u)-h_0)+\frac z2G(u)-\frac z2G(-u) +z(K(u)-k_0),\\
		uK(u)&=z(H(u)-h_0)+z(K(u)-k_0).
\end{align*}
Again we start with the solution of the middle (simplest) function:
\begin{equation*}
	G(u)=\frac {{z}^{2}{u}^{2}f_1+{z}^{2}{u}^{2}e_1+{z}^{2}u-u+{z}^{
			2}uh_0+{z}^{2}uk_0+2 z-{z}^{3}}{{u}^{2}z-u-{z}^{2}u+2 z-{z
		}^{3}}.
\end{equation*}
The denominator is of interest, as one of the factor is `bad':
\begin{equation*}
	u^2z-u-z^2u+2z-z^3=z(u-r_1)(u-r_2),
\end{equation*}
with
\begin{equation*}
	r_1=\frac{1+z^2+\sqrt{1-6z^2+5z^4}}{2z},\quad r_2=\frac{1+z^2-\sqrt{1-6z^2+5z^4}}{2z}.
\end{equation*}
Dividing out the factor $u-r_2$,
\begin{equation*}
	G(u)=\frac {{ r_2} {z}^{2}{ f_1}+{ r_2} {z}^{2}{ e_1}+{z}^
		{2}-1+{z}^{2}{ h_0}+{z}^{2}{ k_0}+u{z}^{2}{ f_0}+u{z}^{2}{ 
			e_1}}{{ r_2} z-1-{z}^{2}+zu},
\end{equation*}
which is what we need.
Dividing out $(u-r_2)(u-2z)u^2$ from the first two functions,
\begin{equation*}
	E(u)=-{\frac { ( { f_1}+{ e_1} ) {z}^{3} ( u+2 z
			) }{{r_2} z+{z}^{3}{r_2}+{z}^{2}{u}^{2}-1-2 {z}^{4}}}
\end{equation*}
and
\begin{equation*}
	F(u)=-{\frac {{z}^{3} ( u{ f_1}+u{ e_1}+1+{ f_1} z+{ e_1} z
			) }{{ r_2} z+{z}^{3}{ r_2}+{z}^{2}{u}^{2}-1-2 {z}^{4}}}.
\end{equation*}
Dividing out $(u-r_2)(u-2z)u$ from the last two functions,
\begin{equation*}
	H(u)={\frac {{z}^{2} ( -uz{ k_0}-uz{ h_0}+{z}^{3}{ e_1}+{z}^{3}{
				f_1}-{z}^{2}{ h_0}-{z}^{2}{ k_0}+{z}^{2}-1+{ f_1} z+{ e_1}
			z ) }{{r_2} z+{z}^{3}{r_2}+{z}^{2}{u}^{2}-1-2 {z}^{4}}}
\end{equation*}
and
\begin{equation*}
	K(u)=-{\frac {{z}^{3} ( { k_0}+{ h_0} )  ( u+2 z
			) }{{r_2} z+{z}^{3}{r_2}+{z}^{2}{u}^{2}-1-2 {z}^{4}}}.
\end{equation*}
Now we can insert $u=0$ and solve:
\begin{align*}
	e_1&=-2 {\frac {{z}^{5}-3 {z}^{3}+2 r_2 {z}^{2}+2 z-r_2}{
			( 1+{z}^{2} )  ( -2+{z}^{2} ) ^{2}}},\\
	f_1&={\frac {-4 z+2 r_2+r_2 {z}^{4}+2 {z}^{3}}{ ( 1+{z}^
			{2} )  ( -2+{z}^{2} ) ^{2}}}
	,\\
	h_0&=-{\frac { ( {z}^{3}-2 z+2 r_2 ) z}{ ( -2+{z}^{2
			} )  ( 1+{z}^{2} ) }}
	,\\
	k_0&=2 {\frac { ( -{z}^{3}+r_2 {z}^{2}+2 z-r_2 ) z}
		{ ( -2+{z}^{2} )  ( 1+{z}^{2} ) }}
	.
\end{align*}
This leads eventually to
\begin{align*}
	E(u)&={\frac { \left( r_2-2 z \right) {z}^{3} \left( u+2 z \right) }
		{ \left( 1+{z}^{2} \right)  (1+2 {z}^{4} -r_2 z-{z}^{3}r_2-{z}^{2}{u}^{2} ) }},\\
	F(u)&=-{\frac { \left( r_2-2 z \right)  \left( r_2 z+{z}^{3}{\it 
				rw}-{z}^{3}u-1 \right) }{ \left( 1+{z}^{2} \right)  (1+2 {z}^{4} -r_2 z-{z}^{3}r_2-{z}^{2}{u}^{2} )}},\\
	G(u)&=-{\frac { \left( r_2-2 z \right)  \left( r_2 {z}^{5}-2 {z
			}^{4}-{z}^{5}u-r_2 z-{z}^{2}+1 \right) }{{z}^{3} \left( 1+{z}^{2
			} \right)  \left( r_2 z-1-{z}^{2}+zu \right) }},\\
	H(u)&=-{\frac { \left( r_2 z+{z}^{3}r_2-{z}^{3}u-1 \right) z
			\left( -3 z+2 r_2 \right) }{ \left( 1+{z}^{2} \right)  (1+2 {z}^{4} -r_2 z-{z}^{3}r_2-{z}^{2}{u}^{2} )}},\\
	K(u)&={\frac { \left( -3 z+2 r_2 \right)  \left( u+2 z \right) {z}^
			{4}}{ \left( 1+{z}^{2} \right) (1+2 {z}^{4} -r_2 z-{z}^{3}r_2-{z}^{2}{u}^{2} ) }}
\end{align*}
and, which is useful for expansion,
\begin{align*}
	E(u)&=\frac{1}{(2-z^2)^2}{\frac { \left( r_2-2 z \right) {z} \left( u+2 z \right) }
		{ \left( 1+{z}^{2} \right)  \bigl(1-\frac{z^2r_2u^2}{(2-z^2)^2}\bigr) }},\\
	F(u)&=\frac{1}{z^2(2-z^2)^2} \frac { \left( r_2-2 z \right)  \left(1- r_2 z-{z}^{3}r_2+{z}^{3}u \right) }{ \left( 1+{z}^{2} \right)  (1-\frac{z^2r_2u^2}{(2-z^2)^2})},\\
	G(u)&={\frac { \left( r_2-2 z \right)  \left( r_2 {z}^{5}-2 {z
			}^{4}-{z}^{5}u-r_2 z-{z}^{2}+1 \right) }{{z}^{3} \left( 1+{z}^{2
			} \right)  \left(1- r_2 z+{z}^{2}-zu \right) }},\\
	H(u)&=\frac{1}{z(2-z^2)^2}{\frac { \left(1 -r_2 z-{z}^{3}r_2+{z}^{3}u \right) 
			\left( -3 z+2 r_2 \right) }{ \left( 1+{z}^{2} \right)  \bigl(1-\frac{z^2r_2u^2}{(2-z^2)^2}\bigr) }},\\
	K(u)&=\frac{1}{(2-z^2)^2}{\frac { \left( -3 z+2 r_2 \right)  \left( u+2 z \right) {z}^
			{2}}{ \left( 1+{z}^{2} \right) \bigl(1-\frac{z^2r_2u^2}{(2-z^2)^2}\bigr)  }}.
\end{align*}
For the simplest function $G(u)$ we show the decomposition in detail:
\begin{align*}
	G(u)&=1+\frac{zu(1+3z^2-zr_2)}{(1+z^2)(1+z^2-zr_2-zu)}\\
	&=1+\frac{u(1+3z^2-zr_2)r_2}{(2-z^2)(1+z^2)\Bigl(1-\frac{r_2}{2-z^2}u\Bigr)}\\
	&=1+\sum_{k\ge1}\Big(\frac{r_2}{2-z^2}\Big)^ku^k \frac{(1+3z^2-zr_2)}{(1+z^2)}.
\end{align*}

Again we can compute the total generating function of paths, regardless where they end:
\begin{align*}
	E(1)+F(1)&+G(1)+H(1)+K(1)\\*
	&=\frac{z(2z^3+9z^2+4z-7)+(z+2)(2z+1)\sqrt{1-6z^2+5z^4}}{2(1+z^2)(1-2z-z^2)}\\
	&=1+z+2z^2+3z^3+5z^4+10z^5+20z^6+38z^7+75z^8+150z^9+ \dots\,.
\end{align*}

\subsection*{The coefficients of $r_2$.}

We report here briefly how to compute the coefficients of 
\begin{equation*}
	r_2^{\ell}=z^\ell\Big(\frac {r_2}{z}\Big)^\ell=z^\ell\bigg(\frac{1+Z-\sqrt{1-6Z+5Z^2}}{2Z}\bigg)^\ell,
\end{equation*}
with $Z=z^2$. As in \cite{garden}, we set $Z=\dfrac{v}{1+3v+v^2}$, so that $\frac{r_2}{z}=2+v$. Then
\begin{align*}
	[Z^n]r_2^{\ell}&=\frac1{2\pi i}\oint \frac{dZ}{Z^{n+1}}(2+v)^\ell\\
	&=\frac1{2\pi i}\oint \frac{dv(1-v^2)}{(1+3v+v^2)^{2}}\frac{(1+3v+v^2)^{n+1}}{v^{n+1}}(2+v)^\ell\\
	&=[v^n](1-v^2)(1+3v+v^2)^{N-1}\sum_{k=0}^\ell\binom{\ell}{k}2^{\ell-k}v^k\\
	&=\sum_{k=0}^\ell\binom{\ell}{k}2^{\ell-k}[v^{n-k}](1-v^2)(1+3v+v^2)^{n-1}\\
	&=\sum_{k=0}^\ell\binom{\ell}{k}2^{\ell-k}\biggl[\binom{n-1;1,3,1}{n-k}-\binom{n-1;1,3,1}{n-k-2}\biggr],
\end{align*}
with weighted \emph{trinomial coefficients} $\dbinom{n;1,3,1}{k}:=[t^k](1+3t+t^2)^n$.

\clearpage

\bibliographystyle{plain}


\end{document}